# THE TRACY–WIDOM LIMIT FOR THE LARGEST EIGENVALUES OF SINGULAR COMPLEX WISHART MATRICES

By Alexei Onatski

*Columbia University*

This paper extends the work of El Karoui [*Ann. Probab.* **35** (2007) 663–714] which finds the Tracy–Widom limit for the largest eigenvalue of a nonsingular $p$-dimensional complex Wishart matrix $W_{\mathbb{C}}(\Omega_p, n)$ to the case of several of the largest eigenvalues of the possibly singular $(n < p)$ matrix $W_{\mathbb{C}}(\Omega_p, n)$. As a byproduct, we extend all results of Baik, Ben Arous and Peche [*Ann. Probab.* **33** (2005) 1643–1697] to the singular Wishart matrix case. We apply our findings to obtain a 95% confidence set for the number of common risk factors in excess stock returns.

**1. Introduction.** The goal of this paper is to establish the joint asymptotic distribution of a finite number of properly scaled and centered largest eigenvalues of a $p$-dimensional complex Wishart matrix $W_{\mathbb{C}}(\Omega_p, n)$ as both $n$ and $p$ tend to infinity in such a way that $n/p$ remains in a compact subset of $(0, \infty)$. The paper extends Baik, Ben Arous and Peche [4] and El Karoui [9], which find the asymptotic distribution of the scaled and centered single largest eigenvalue of $W_{\mathbb{C}}(\Omega_p, n)$ under the assumption that $n/p$ remains in a compact subset of $[1, \infty)$. When $n/p$ is less than 1, the Wishart matrix is singular. The main contribution of this paper is extending [4] and [9] to the singular Wishart case.

One need for such an extension arises in a companion paper [17], which develops statistical tests of various hypotheses concerning the number of factors in Chamberlain and Rothschild's [8] approximate factor model, which is widely used in empirical macroeconomics and finance. The model considers a double-infinite sequence of random variables $\{\xi_{it}, i, t \in \mathbb{N}\}$ such that, for any $i, t \in \mathbb{N}$,

$$\xi_{it} = \Lambda'_i F_t + \eta_{it}, \tag{1}$$

where $F_t$ and $\Lambda_i$ are $k$-dimensional ($k < \infty$) vectors of unobserved common factors and factor loadings, respectively, and $\eta_{it}$ is an unobserved idiosyncratic component of $\xi_{it}$. In contrast to the classical factor model (see [1], Chapter 14), the idiosyncratic components are allowed to be correlated over the $i$th dimension. The identification of the idiosyncratic components is achieved by assuming that the largest eigenvalue of the covariance matrix of vector $\{\eta_{it}\}_{1 \leq i \leq p}$ stays bounded as $p$ tends to infinity, whereas the smallest nonzero eigenvalue of the covariance matrix of vector $\{\Lambda_i' F_t\}_{1 \leq i \leq p}$ diverges to infinity as fast as $p$. These assumptions are often interpreted as formalizing the requirements that the common factors nontrivially influence all data points, whereas the idiosyncratic components have only local effects.

Researchers in macroeconomics and finance use Chamberlain and Rothschild's [8] model to handle high-dimensional data sets. They interpret $F_t$ as a vector of factors nontrivially influencing hundreds of available macroeconomic indicators (see, e.g., [25]) or, in the case of finance, as a vector of the risk factors common to hundreds of stock returns (see, e.g., [7]). An important practical question is how many such factors there are. Under the assumptions of the Chamberlain–Rothschild model, one can equivalently ask how many eigenvalues of $XX'/n$, where $X$ is the data matrix $\{\xi_{it}\}_{1 \leq i \leq p, 1 \leq t \leq n}$, diverge to infinity as both $p$ and $n$ tend to infinity in such a way that $n/p$ remains in a compact subset of $(0, \infty)$.

Using eigenvalue perturbation theory, [17] shows that if the true number of factors is $k_0$, then the asymptotic distribution of the scaled and centered $(k_0 + 1)$th, $(k_0 + 2)$th, and so on, eigenvalues of $XX'/n$ is the same as that of the scaled and centered 1st, 2nd, and so on, eigenvalues of the sample covariance matrix of the idiosyncratic components (slightly abusing a standard definition, we define the sample covariance matrix of vectors $v_1, \ldots, v_n$ as $\sum_{i=1}^{n} v_i v_i'/n$). Therefore, assuming that $\{\eta_{it}\}_{1 \leq i \leq p}$ are i.i.d. (over $t \in \mathbb{N}$) complex Gaussian $N_{\mathbb{C}}(0, \Sigma_p)$ vectors (to make such an assumption realistic, we perform a preliminary transformation of real-valued data into a complex-valued form), a test (described in more detail in Section 4 below) of the null hypothesis that the true number of factors equals $k_0$, against an alternative of more than $k_0$ factors, can be based on checking whether the $(k_0 + 1)$th, $(k_0 + 2)$th, and so on, eigenvalues of the sample covariance matrix of the data are drawn from the joint distribution of the largest eigenvalues of $W_{\mathbb{C}}(\Sigma_p/n, n)$. Since, in macroeconomics and finance, the cross-sectional dimension $p$ of data is often larger than their time-series dimension $n$, to obtain the asymptotic critical values of the test, we must analyze the joint asymptotic distribution of the largest eigenvalues of a singular complex Wishart matrix.

El Karoui [9] proves that the asymptotic distribution of the properly scaled and centered largest eigenvalue of a nonsingular complex Wishart matrix, $W_\mathbb{C}(\Sigma_p/n, n)$, is the Tracy–Widom distribution of type two (TW$_2$). TW$_2$ refers to a distribution with the cumulative distribution function

$$F(x) \equiv \exp\left(-\int_x^\infty (x-s)q^2(s)\,ds\right),$$

where $q(s)$ is the solution of an ordinary differential equation

$$q''(s) = sq(s) + 2q^3(s),$$

which is asymptotically equivalent to the Airy function $Ai(s)$ (see [15]) as $s \to \infty$. It plays an important role in large random matrix theory (see [14]) because it is the asymptotic distribution of the scaled and centered largest eigenvalue of a matrix from the so-called Gaussian Unitary Ensemble (GUE) as the size of the matrix tends to infinity.

The GUE is the collection of all $N \times N$ Hermitian matrices with i.i.d. complex Gaussian $N_\mathbb{C}(0, 1/N)$ lower triangular entries and (independent of them) i.i.d. real Gaussian $N(0, 1/N)$ diagonal entries. Let $d_1 \geq \cdots \geq d_N$ be eigenvalues of a matrix from the GUE. Define $\tilde{d}_i = N^{2/3}(d_i - 2)$. Tracy and Widom [26] studied the asymptotic distribution of a few of the largest eigenvalues of matrices from the GUE when $N \to \infty$. They described the asymptotic marginal distributions of $\tilde{d}_i$, $i = 1, \ldots, m$, where $m$ is any fixed positive integer, in terms of a solution of a completely integrable system of partial differential equations. If we are interested in the asymptotic distribution of the largest eigenvalue only, then the system simplifies to the single ordinary differential equation given above.

In this paper, we extend El Karoui's [9] results to show that the asymptotic distribution of the scaled and centered $m$ largest eigenvalues ($m < \infty$) of a possibly singular complex Wishart matrix $W_\mathbb{C}(\Sigma_p/n, n)$ is the same as the joint asymptotic distribution of $\tilde{d}_1, \ldots, \tilde{d}_m$. We follow [24] in calling such a joint distribution the *joint Tracy–Widom distribution*.

Random matrix theory has developed the following powerful method of analysis of the joint asymptotic distribution of a few of the largest eigenvalues of various random matrices as the dimensionality of the matrices tends to infinity. First, the joint distribution of a few of the largest eigenvalues is expressed in terms of the probabilities $P(i_1, \ldots, i_m; J_1, \ldots, J_m)$ that disjoint subsets $J_1, \ldots, J_m$ of the real line contain exactly $i_1, \ldots, i_m$ eigenvalues. The latter probabilities are then represented in the form of Fredholm determinants of operators indexed by the dimensionality of the random matrix under consideration. Finally, it is proved that the operators converge in the trace-class norm as the dimensionality tends to infinity and the corresponding limits are found. The outcome of such an analysis is an expression of the joint distribution of a few of the largest eigenvalues in terms of Fredholm determinants of the limiting integral operators. Often, kernels of these operators have a relatively simple form which ensures further detailed analysis of the joint distribution of the largest eigenvalues.

Let $\lambda_j$ be the $j$th largest eigenvalue of a complex Wishart matrix $W_{\mathbb{C}}(\Sigma_p/n,n)$. The first step of the above method is then performed as follows. For any real $s_1 > \cdots > s_m > 0$,

(2)
$$\Pr(\lambda_1 \leq s_1, \ldots, \lambda_m \leq s_m)$$
$$= \sum_I P(i_1, \ldots, i_m; (s_1, \infty), (s_2, s_1], \ldots, (s_m, s_{m-1}]),$$

where $I$ consists of all sets of $m$ nonnegative integers $i_1, \ldots, i_m$ such that $i_1 = 0$ and $i_{j+1} \leq j - i_j - \cdots - i_1$ for $j = 1, \ldots, m-1$. In the special case where only the largest eigenvalue is analyzed, we have $\Pr(\lambda_1 \leq s_1) = P(0, (s_1, \infty)) = E \prod_{j=1}^p [1 - \chi_{(s_1, \infty)}(\lambda_j)]$, where $\chi_J(\lambda)$ denotes the indicator function of the set $J$ and the expectation is with respect to the joint distribution of $\lambda_1, \ldots, \lambda_p$. For this special case, Baik, Ben Arous and Peche [4] perform the second step of the above method. Assuming that $n \geq p$, they show that $E \prod_{j=1}^p [1 - \chi_{(s_1, \infty)}(\lambda_j)]$ equals the Fredholm determinant $\det(1 - K_{n,p})$, where $K_{n,p}$ is an operator acting on $L^2((s_1, \infty))$ with a kernel that has a convenient integral representation.

To establish the asymptotic distribution of the largest eigenvalue of $W_{\mathbb{C}}(\Sigma_p/n,n)$, El Karoui [9] starts from Baik, Ben Arous and Peche [4] result. He then finds centering and scaling sequences $\mu_{n,p}$ and $\sigma_{n,p}$ such that, as both $n$ and $p$ tend to infinity, the recentered and rescaled version of the operator $K_{n,p}$, $S_{n,p}$, converges in the trace-class norm to an operator $E \cdot Ai \cdot E^{-1}$ acting on $L^2((s_1, \infty))$, where $E$ is an operator which entails multiplication by a certain function and $Ai$ is an integral operator with the Airy kernel $Ai(x,y) = \int Ai(x+u)Ai(y+u)\,du$, where $Ai(x)$ is the Airy function. Since the Fredholm determinant is continuous with respect to the trace-class norm and since it is invariant with respect to conjugation, El Karoui [9] concludes that the distribution of the centered and scaled largest eigenvalue of $W_{\mathbb{C}}(\Sigma_p/n,n)$ converges to the distribution defined by $\det(I - Ai)$, which is $TW_2$ (see [26]).

A careful inspection of El Karoui's [9] proofs reveals that they only involve the assumption that $n/p$ remains in a compact subset of $[1, \infty)$ to be able to use the determinantal representation of the cumulative distribution function of the largest eigenvalue of $W_{\mathbb{C}}(\Sigma_p/n,n)$ established by Baik, Ben Arous and Peche [4]. Therefore, we first extend Baik, Ben Arous and Peche [4] to the case of a singular complex Wishart matrix. Somewhat unexpectedly, we find that not only the determinantal representation of the cumulative distribution function of the largest eigenvalue of $W_{\mathbb{C}}(\Sigma_p/n,n)$, but also all the rest of their results, hold for the singular Wishart case without any extra qualifications. Our extension of [9] easily follows from the extension of [4].

The rest of this paper is organized as follows. In Section 2, we prove our generalization of [4]. Section 3 generalizes [9] to the case of several eigenvalues of a possibly singular complex Wishart matrix. Section 4 contains an application to the determination of the number of common risk factors in stock return data. Section 5 concludes. The Appendix contains proofs of some of the tangential statements of this paper.

**2. Extension of Baik, Ben Arous and Peche [4].** In this section, we extend Baik, Ben Arous and Peche [4] analysis to the singular situation when $n < p$ and to the case of several of the largest eigenvalues of $W_{\mathbb{C}}(\Sigma_p/n, n)$. Note that for a general positive integer $m$, we have

$$
(3) \quad \begin{aligned} & P(i_1, \ldots, i_m; J_1, \ldots, J_m) \\ & = \frac{1}{i_1! \cdots i_m!} \frac{\partial^{i_1 + \cdots + i_m}}{\partial z_1^{i_1} \cdots \partial z_m^{i_m}} E \prod_{j=1}^{p} \left[ 1 + \sum_{k=1}^{m} (z_k - 1) \chi_{J_k}(\lambda_j) \right] \Bigg|_{z_1 = \cdots = z_m = 0} \end{aligned}
$$

(see, e.g., formula (4.1) in [27]). Below, we establish a determinantal representation of $E \prod_{j=1}^{p}[1 + \sum_{k=1}^{m}(z_k - 1)\chi_{J_k}(\lambda_j)]$ which does not depend on whether $n < p$ or $n \geq p$.

First, for the case $n < p$, we will find a convenient expression for the joint density of the nonzero eigenvalues $\lambda_1, \ldots, \lambda_n$ of a singular complex Wishart matrix $W_{\mathbb{C}}(\Sigma_p/n, n)$. Let $\pi_j$ be the *inverse* of the $j$th largest eigenvalue of $\Sigma_p$. Define $\vec{\lambda} = (\lambda_1, \ldots, \lambda_n)'$, $\vec{\pi} = (\pi_1, \ldots, \pi_p)'$, $\Lambda = \mathrm{diag}(\lambda_1, \ldots, \lambda_n)$ and $\Pi = \mathrm{diag}(\pi_1, \ldots, \pi_p)$. As shown in [18], formula (25), the joint density equals

$$
(4) \quad f(\vec{\lambda}) = \mathrm{const} \cdot V(\vec{\lambda})^2 \prod_{j=1}^{n} \lambda_j^{p-n} \int_{Q_1 \in CV(p,n)} e^{-n \, \mathrm{tr}(\Pi Q_1 \Lambda Q_1^*)} (Q^* \, dQ_1),
$$

where $V(\vec{\lambda}) = \prod_{1 \leq i < j \leq n}(\lambda_j - \lambda_i)$, $CV(p, n)$ denotes the complex Stiefel manifold of $p \times n$ matrices with orthonormal columns and $(Q^* \, dQ_1)$ is the exterior differential form representing the uniform measure on the complex Stiefel manifold. Throughout the remainder of this paper, "const" denotes possibly different constants that may depend on $p$, $n$ and $\vec{\pi}$, but not on $\vec{\lambda}$.

Note that

$$
(5) \quad \int_{Q_1 \in CV(p,n)} e^{-n \, \mathrm{tr}(\Pi Q_1 \Lambda Q_1^*)} (Q^* \, dQ_1) = \frac{\int_{R \in U(p)} e^{-n \, \mathrm{tr}(\Pi R_1 \Lambda R_1^*)} (R^* \, dR)}{\mathrm{Vol}\{U(p - n)\}},
$$

where $R_1$ is a $p \times n$ matrix of the first $n$ columns of the matrix $R \equiv [R_1, R_2]$, and $U(p)$ is the set of all $p \times p$ unitary matrices. A proof of (5) can be found in the Appendix. Now, since the unitary group is compact, we have

$$
(6) \quad \int_{R \in U(p)} e^{-n \, \mathrm{tr}(\Pi R_1 \Lambda R_1^*)} (R^* \, dR) = \lim_{\varepsilon \to 0} \int_{R \in U(p)} e^{-n \, \mathrm{tr}(\Pi R \Lambda_\varepsilon R^*)} (R^* \, dR),
$$

where $\Lambda_\varepsilon := \mathrm{diag}(\lambda_1, \ldots, \lambda_n, \varepsilon, 2\varepsilon, \ldots, \varepsilon(p - n))$.

The integral on the right-hand side of (6) is called the Harish–Chandra–Itzykson–Zuber integral (see [14], Appendix A5 and page 648). It can be simplified as follows. Define $\vec{\lambda}_\varepsilon = (\lambda_1, \ldots, \lambda_n, \varepsilon, 2\varepsilon, \ldots, \varepsilon(p-n))'$ and let $\lambda_{\varepsilon k}$ be the $k$th component of vector $\vec{\lambda}_\varepsilon$. We then have

$$
\begin{aligned}
(7) \quad \int_{R \in U(p)} & e^{-n \operatorname{tr}(\Pi R \Lambda_\varepsilon R^*)} (R^* \, dR) \\
&= \operatorname{const} \cdot (V(\vec{\pi}) V(\vec{\lambda}_\varepsilon))^{-1} \det(e^{-n \pi_j \lambda_{\varepsilon k}})_{1 \leq j,k \leq p}.
\end{aligned}
$$

Here, we assume that all of the $\pi_i$ are different. If some $\pi_i$ are equal, then the formula should be changed according to l'Hôpital's theorem.

Let $\{i_1, \ldots, i_p\}$ be a set of indices equal to the set $\{1, 2, \ldots, p\}$ and such that $i_1 < \cdots < i_n$ and $i_{n+1} < \cdots < i_p$. Denote the multi-index $(i_1, \ldots, i_n)$ as $\alpha$ and the multi-index $(i_{n+1}, \ldots, i_p)$ as $\bar{\alpha}$, and let $x_\alpha$ denote $(x_{i_1}, \ldots, x_{i_n})'$, $x_{\alpha(k)}$ denote $x_{i_k}$, $x_{\bar{\alpha}}$ denote $(x_{i_{n+1}}, \ldots, x_{i_p})'$ and $x_{\bar{\alpha}(k)}$ denote $x_{i_{n+k}}$. Finally, let $|\alpha|$ denote $i_1 + \cdots + i_n$. Then, by the Laplace expansion theorem, $\det(e^{-n\pi_j \lambda_{\varepsilon k}})_{1 \leq j,k \leq p}$ is equal in absolute value to

$$
(8) \quad \sum_\alpha (-1)^{|\alpha|} \det(e^{-n \pi_{\alpha(j)} \lambda_k})_{1 \leq j,k \leq n} \det((e^{-n \pi_{\bar{\alpha}(j)} \varepsilon})^k)_{1 \leq j,k \leq p-n}.
$$

The second determinant in the above sum is a Vandermonde determinant. Hence,

$$
\begin{aligned}
(9) \quad \det&((e^{-n \pi_{\bar{\alpha}(j)} \varepsilon})^k)_{1 \leq j,k \leq p-n} \\
&= e^{-n\varepsilon \sum_{j=1}^{p-n} \pi_{\bar{\alpha}(j)}} \prod_{1 \leq j < k \leq p-n} (e^{-n \pi_{\bar{\alpha}(k)} \varepsilon} - e^{-n \pi_{\bar{\alpha}(j)} \varepsilon}).
\end{aligned}
$$

Further, note that

$$
(10) \quad V(\vec{\lambda}_\varepsilon) = \operatorname{const} \cdot \varepsilon^{\binom{p-n}{2}} V(\vec{\lambda}) \prod_{i=1}^n \prod_{k=n+1}^p (\varepsilon(k-n) - \lambda_i).
$$

Combining (5)–(10) and taking the limit as $\varepsilon \to 0$, we obtain

$$
\begin{aligned}
(11) \quad \int_{Q_1 \in CV(p,n)} & e^{-n \operatorname{tr}(\Pi Q_1 \Lambda Q_1^*)} (Q^* \, dQ_1) \\
&= \operatorname{const} \cdot (V(\vec{\pi}) V(\vec{\lambda}))^{-1} \\
&\quad \times \prod_{i=1}^n \lambda_i^{n-p} \sum_\alpha (-1)^{|\alpha|} V(\pi_{\bar{\alpha}}) \det(e^{-n \pi_{\alpha(j)} \lambda_k})_{1 \leq j,k \leq n}.
\end{aligned}
$$

Substituting (11) into (4), we find that the joint density of the nonzero eigenvalues $\vec{\lambda} = (\lambda_1, \ldots, \lambda_n)$ of the singular complex Wishart matrix $W_{\mathbb{C}}(\Sigma_p/n, n)$ equals

$$
(12) \quad f(\vec{\lambda}) = \operatorname{const} \cdot V(\vec{\lambda}) \sum_\alpha (-1)^{|\alpha|} V(\pi_{\bar{\alpha}}) \det(e^{-n \pi_{\alpha(j)} \lambda_k})_{1 \leq j,k \leq n}.
$$

We are now ready to generalize Proposition 2.1 of [4], which establishes a determinantal representation of $E\prod_{j=1}^{p}[1 - \chi_{(s_1,\infty)}(\lambda_j)]$, to the case of $E\prod_{j=1}^{p}[1 + \sum_{k=1}^{m}(z_k - 1)\chi_{J_k}(\lambda_j)]$ and general $n$ and $p$.

PROPOSITION 1. *For any fixed $q$ satisfying $0 < q < \min\{\pi_j\}_{j=1}^{p}$, let $K_{n,p}$ be the operator acting on $L^2((0,\infty))$ with kernel*

$$K_{n,p}(\eta, \zeta) = \frac{n}{(2\pi i)^2} \int_\Gamma dz \int_\Sigma dw \, e^{-\eta n(z-q) + \zeta n(w-q)} \frac{1}{w-z} \left(\frac{z}{w}\right)^n \prod_{k=1}^{p} \frac{\pi_k - w}{\pi_k - z},$$

*where $\Sigma$ is a simple closed contour enclosing $0$ and lying in $\{w : \mathrm{Re}(w) < q\}$ and $\Gamma$ is a simple closed contour enclosing $\pi_1, \ldots, \pi_p$ and lying in $\{z : \mathrm{Re}(z) > q\}$, both oriented counterclockwise. Then, for any real-valued measurable bounded function $f(x)$ which equals $0$ for any $x \leq 0$,*

$$E\prod_{j=1}^{p}(1 + f(\lambda_j)) = \det(1 + K_{n,p}f),$$

*where $1$ is the identity operator and $f$ is the operator which entails multiplication by the function $f(\cdot)$.*

PROOF. Let us first focus on the case when $n < p$. The eigenvalues $\lambda_{n+1}, \ldots, \lambda_p$ then equal zero and we have $E\prod_{j=1}^{p}(1 + f(\lambda_j)) = E\prod_{j=1}^{n}(1 + f(\lambda_j))$. Using the equality $V(\vec{\lambda}) = \det(\lambda_k^{j-1})_{1 \leq j,k \leq n}$ and (12), we obtain

$$E\prod_{j=1}^{n}(1 + f(\lambda_j)) = \mathrm{const} \cdot \sum_\alpha (-1)^{|\alpha|} V(\pi_{\bar{\alpha}})$$
$$\times \int_0^\infty \cdots \int_0^\infty \det(\lambda_k^{j-1})$$
$$\times \det(e^{-n\pi_{\alpha(j)}\lambda_k}) \prod_{k=1}^{n}(1 + f(\lambda_k)) \, d\lambda_k.$$

Using Andreief's [2] identity $\int \cdots \int \det(f_j(x_k)) \det(g_j(x_k)) \prod_k d\mu(x_k) = \det(\int f_j(x) g_k(x) \, d\mu(x))$, we find that

$$E\prod_{j=1}^{n}(1 + f(\lambda_j))$$
(13) $$= \mathrm{const} \cdot \sum_\alpha (-1)^{|\alpha|} V(\pi_{\bar{\alpha}})$$
$$\times \det\left(\int_0^\infty (1 + f(\lambda))\lambda^{j-1} e^{-n\pi_{\alpha(k)}\lambda} \, d\lambda\right)_{1 \leq j,k \leq n}.$$

Now, define $\nu = n - p$ (note that it is less than zero) and set

$$\phi_j(\lambda) = \begin{cases} 0, & \text{if } j \leq -\nu, \\ \dfrac{n^{j+\nu}}{\Gamma(j+\nu)} \lambda^{j-1+\nu} e^{-nq\lambda}, & \text{if } j > -\nu, \end{cases} \quad \text{for } j = 1, \ldots, p,$$

$$\Phi_k(\lambda) = e^{-n(\pi_k - q)\lambda}, \qquad \text{for } k = 1, \ldots, p,$$

for any $0 < q < \min\{\pi_j\}_{j=1}^p$. Also, let

$$A = (A_{jk})_{1 \leq j,k \leq p}, \qquad A_{jk} = \pi_k^{-j-\nu}.$$

Note that $A_{jk} = \int_0^\infty \phi_j(\lambda) \Phi_k(\lambda) \, d\lambda$ for $j > -\nu$. Since $A$ is a simple modification of a Vandermonde matrix, we have

$$(14) \qquad \det A = \prod_{j=1}^p \frac{1}{\pi_j^{\nu+1}} \prod_{1 \leq j < k \leq p} (\pi_k^{-1} - \pi_j^{-1}).$$

Thus, $A$ is invertible when all of the $\pi_i$ are distinct.

Next, define operators $B : L^2((0, \infty)) \to l^2(\{1, \ldots, p\})$ and $C : l^2(\{1, \ldots, p\}) \to L^2((0, \infty))$ by

$$B(j, \lambda) = \phi_j(\lambda), \qquad C(\lambda, k) = \Phi_k(\lambda)$$

and let $f$ be the operator $L^2((0, \infty)) \to L^2((0, \infty))$ which entails multiplication by a real-valued measurable bounded function $f(x)$. Then, since

$$(15) \qquad \frac{n^j}{\Gamma(j)} \int_0^\infty f(\lambda) \lambda^{j-1} e^{-n\pi_{\alpha(k)}\lambda} \, d\lambda = (BfC)(j - \nu, \alpha(k))$$

for $j = 1, \ldots, n$, we find, using (13), that

$$(16) \qquad E \prod_{j=1}^n (1 + f(\lambda_j)) = \text{const} \cdot \sum_\alpha (-1)^{|\alpha|} V(\pi_{\bar{\alpha}}) \det(A^{(\alpha)} + B^{(\alpha)} f C^{(\alpha)}),$$

where $A^{(\alpha)}$ is a submatrix of $A$ that consists of the intersection of its columns numbered $\alpha(1), \ldots, \alpha(n)$ with its last $n$ rows, that is, $A^{(\alpha)} = (A_{j-\nu,\alpha(k)})_{1 \leq j,k \leq n}$. Similarly, $B^{(\alpha)}$ is an operator with the kernel that consists of the last $n$ elements of the kernel of $B$ and $C^{(\alpha)}$ is an operator with the kernel that consists of elements, numbered $\alpha(1), \ldots, \alpha(n)$, of the kernel of $C$.

Note that the right-hand side of (16) is proportional to the Laplace expansion of $\det(A - BfC)$. To see this, use (15) and observe that the kernel of $A - BfC$ has the following form:

$$\begin{pmatrix} \pi_1^{-1-\nu} & \cdots & \pi_p^{-1-\nu} \\ \vdots & & \vdots \\ 1 & \cdots & 1 \\ \pi_1^{-1} + \frac{n}{\Gamma(1)} \int_0^\infty f(\lambda) e^{-n\pi_1 \lambda} d\lambda & \cdots & \pi_p^{-1} + \frac{n}{\Gamma(1)} \int_0^\infty f(\lambda) e^{-n\pi_p \lambda} d\lambda \\ \vdots & & \vdots \\ \pi_1^{-n} + \frac{n^n}{\Gamma(n)} \int_0^\infty f(\lambda) \lambda^{n-1} e^{-n\pi_1 \lambda} d\lambda & \cdots & \pi_p^{-n} + \frac{n^n}{\Gamma(n)} \int_0^\infty f(\lambda) \lambda^{n-1} e^{-n\pi_p \lambda} d\lambda \end{pmatrix}.$$

Hence,

$$E \prod_{j=1}^n (1 + f(\lambda_j)) = \text{const} \cdot \det(A + BfC)$$
$$= \text{const} \cdot \det(A^{-1}) \det(1 + A^{-1} BfC).$$

So, interchanging the order of the composition of operators under the determinant, we find that

$$E \prod_{j=1}^n (1 + f(\lambda_j)) = \text{const} \cdot \det(1 + CA^{-1} Bf).$$

By setting $f(\cdot)$ equal to minus the indicator function of $(s, \infty)$ and letting $s \to \infty$ in both sides of the above equality, we find that "const" in the above formula equals 1. Thus,

$$E \prod_{j=1}^n (1 + f(\lambda_j)) = \det(1 + CA^{-1} Bf).$$

From this point on the proof of Baik, Ben Arous and Peche [4] progresses practically without changes. We will provide it here in order to make this paper self-contained. The kernel of the operator $CA^{-1}B$ in the above determinant is

$$CA^{-1}B(\eta, \zeta) = \sum_{j=1}^p C(\eta, k)(A^{-1}B)(j, \zeta), \qquad \eta, \zeta > 0.$$

Further, from Cramér's rule,

(17) $$(A^{-1}B)(j, \zeta) = \frac{\det A^{(j)}(\zeta)}{\det A},$$

where $A^{(j)}(\zeta)$ is the matrix given by $A$ with the $j$th column replaced by the vector $(\phi_1(\zeta), \ldots, \phi_p(\zeta))'$. To compute $\det A^{(j)}$, note that

$$\frac{1}{2\pi i} \int_\Sigma \frac{e^w}{w^a} dw = \begin{cases} \dfrac{1}{\Gamma(a)}, & \text{if } a \text{ is positive integer,} \\ 0, & \text{if } a \text{ is zero or negative,} \end{cases}$$

where $\Sigma$ is any simple closed contour with counterclockwise orientation enclosing the origin 0. By replacing $w \to \zeta n w$ and setting $a = k + \nu$, this implies that

$$\frac{\zeta^{-(k-1+\nu)}}{2\pi i} \int_\Sigma e^{\zeta n w} \frac{n}{(nw)^{k+\nu}} dw = \begin{cases} \dfrac{1}{\Gamma(k+\nu)}, & \text{if } k > -\nu, \\ 0, & \text{if } k \leq -\nu, \end{cases}$$

and therefore that

$$\phi_k(\zeta) = \frac{1}{2\pi i} \int_\Sigma e^{\zeta n(w-q)} \frac{n}{w^{k+\nu}} dw.$$

Substituting this formula for $\phi_k(\zeta)$ in the $j$th column of $A^{(j)}$ and pulling out the integrals over $w$, we obtain

$$\det A^{(j)}(\zeta) = \frac{1}{2\pi i} \int_\Sigma e^{\zeta n(w-q)} \det(A'(w)) n \, dw,$$

where the entries of $A'(w)$ are $A'_{ab}(w) = 1/p_b^{a+\nu}$, where $p_b = \pi_b$ when $b \neq j$ and $p_b = w$ when $b = j$. Hence, by the formula for a Vandermonde determinant,

$$\det A^{(j)}(\zeta) = \prod_{k \neq j} \frac{1}{\pi_k^{1+\nu}} \frac{1}{2\pi i} \int_\Sigma e^{\zeta n(w-q)} \prod_{1 \leq a < b \leq n} (p_b^{-1} - p_a^{-1}) \frac{n \, dw}{w^{1+\nu}}$$

and so, using (14) and (17), we obtain

$$(A^{-1}B)(j,\zeta) = \frac{n \pi_j^{p+\nu}}{2\pi i} \int_\Sigma e^{\zeta n(w-q)} \prod_{k \neq j} \frac{w - \pi_k}{\pi_j - \pi_k} \frac{dw}{w^{p+\nu}}.$$

However, for any simple closed contour $\Gamma$ that encloses $\pi_1, \ldots, \pi_p$ but excludes $w$ and which is oriented counterclockwise,

$$\frac{1}{2\pi i} \int_\Gamma z^n e^{-\eta n z} \frac{1}{w-z} \prod_{k=1}^p \frac{w - \pi_k}{z - \pi_k} dz = \sum_{j=1}^p \pi_j^n e^{-n\pi_j \eta} \prod_{k \neq j} \frac{w - \pi_k}{\pi_j - \pi_k}.$$

Therefore, we find

$$CA^{-1}B(\eta, \zeta)$$
$$= \frac{n}{(2\pi i)^2} \int_\Gamma dz \int_\Sigma dw \, e^{-\eta n(z-q) + \zeta n(w-q)} \frac{1}{w-z} \prod_{k=1}^p \frac{w - \pi_k}{z - \pi_k} \left(\frac{z}{w}\right)^n,$$

which completes the proof when all of the $\pi_j$ are distinct. When some $\pi_j$ are equal, the formula for the kernel $CA^{-1}B(\eta,\zeta)$ follows by taking proper limits and using l'Hôpital's theorem.

For the case when $n \geq p$ and $f(x)$ equals minus the indicator function for the interval $(s, \infty)$, where $s \in R$, the proposition is equivalent to Proposition 2.1 of [4]. Extending Baik, Ben Arous and Peche's [4] proof to the case of general $f(x)$ while keeping their assumption that $n \geq p$ is straightforward. To save space, we omit such an extension from the proof. □

In the next section, we will use Proposition 1 to extend the results of El Karoui [9] to the case of several of the largest eigenvalues of a complex singular Wishart matrix. Concluding this section, we would like to note that our extension of Proposition 2.1 of [4] implies that the main results of that paper, namely Theorems 1.1 and 1.2, hold under the assumption that $n/p$ ($M/N$ in the notation of [4]) remains in a compact subset of $(0, +\infty)$ as both $n$ and $p$ tend to infinity. This assumption relaxes Baik, Ben Arous and Peche's [4] requirement that $n/p$ remains in a compact subset of $[1, +\infty)$. The Appendix contains a brief list of changes that should be made to the proofs of Baik, Ben Arous and Peche [4] (beyond the extension of Proposition 2.1) to justify such a relaxation.

**3. Extension of El Karoui [9].** In this section, we will prove that the joint distribution of the first $m$ scaled and centered eigenvalues of a complex Wishart matrix $W_{\mathbb{C}}(\Sigma_p/n, n)$ weakly converges to the $m$-dimensional joint Tracy–Widom distribution. Such a convergence takes place in both cases $n \geq p$ and $n < p$. The scaling and centering sequences are the same for all of the $m$ eigenvalues and have the form proposed in [9].

PROPOSITION 2. *Let $H_p$ be the spectral distribution of $\Sigma_p$. Let $c_{n,p}$ be the unique solution in $[0, \pi_1)$ of the equation*

$$\int \left(\frac{\lambda c_{n,p}}{1 - \lambda c_{n,p}}\right)^2 dH_p(\lambda) = \frac{n}{p}.$$

*Assume that $n/p$ remains in a compact subset of $(0, \infty)$, $\limsup \pi_1^{-1} < \infty$, $\liminf \pi_p^{-1} > 0$ and $\limsup c_{n,p}/\pi_1 < 1$. Define*

$$\mu_{n,p} = \frac{1}{c_{n,p}}\left(1 + \frac{p}{n}\int \frac{\lambda c_{n,p}}{1 - \lambda c_{n,p}} dH_p(\lambda)\right)$$

*and*

$$\sigma_{n,p} = \frac{1}{n^{2/3} c_{n,p}}\left(1 + \frac{p}{n}\int \left(\frac{\lambda c_{n,p}}{1 - \lambda c_{n,p}}\right)^3 dH_p(\lambda)\right)^{1/3}.$$

Then, as $n$ and $p$ tend to infinity, the joint distribution of the first $m$ centered and scaled eigenvalues $\sigma_{n,p}^{-1}(\lambda_1 - \mu_{n,p}), \ldots, \sigma_{n,p}^{-1}(\lambda_m - \mu_{n,p})$ of matrix $W_{\mathbb{C}}(\Sigma_p/n, n)$ weakly converges to the $m$-dimensional joint Tracy–Widom distribution.

PROOF. For a short proof of the uniqueness of $c_{n,p}$, see [9], formula (11) and the discussion that follows the formula. Let us first prove that $\liminf c_{n,p} > 0$ and $\limsup c_{n,p} < \infty$. Since, by assumption, $\liminf n/p > 0$, there exists $\gamma > 0$ such that $n/p > \gamma^2$. We have

$$\gamma^2 < \frac{n}{p} = \int \left( \frac{\lambda c_{n,p}}{1 - \lambda c_{n,p}} \right)^2 dH_p(\lambda) \leq \left( \frac{c_{n,p}/\pi_1}{1 - c_{n,p}/\pi_1} \right)^2$$

and, therefore, $c_{n,p} > \pi_1 \frac{\gamma}{1+\gamma}$. This implies that $\liminf c_{n,p} > 0$ because, by assumption, $\limsup \pi_1^{-1} < \infty$. Further, since $c_{n,p} < \pi_1 \leq \pi_p$, the assumption that $\liminf \pi_p^{-1} > 0$ implies that $\limsup c_{n,p} < \infty$. Note that the facts just established, that $\liminf c_{n,p} > 0$ and $\limsup c_{n,p} < \infty$, and the assumptions that $\limsup c_{n,p}/\pi_1 < 1$ and that $n/p$ remains in a compact subset of $(0, \infty)$ imply that $\mu_{n,p}$ remains in a compact subset of $(0, \infty)$ and that $\sigma_{n,p}$ decays to zero as fast as $n^{-2/3}$ when $n$ tends to infinity.

Now, let $x_1 > \cdots > x_m$ be any real numbers. Since $\mu_{n,p}$ remains in a compact subset of $(0, \infty)$, whereas $\sigma_{n,p}$ tends to zero as $n \to \infty$, there exists $N > 0$ such that for any $n > N$, $s_i = \mu_{n,p} + \sigma_{n,p} x_i$ $(i = 1, \ldots, m)$ are positive numbers. In what follows, we will always take $n > N$. Consider the function $f_{z_1,\ldots,z_m}(x) = \sum_{k=1}^{m}(z_k - 1)\chi_{J_k}(x)$, where $J_1, \ldots, J_m$ equal $(s_1, \infty), (s_2, s_1], \ldots, (s_m, s_{m-1}]$, respectively, and $z_1, \ldots, z_m$ are any complex numbers. According to (2) and (3),

$$\Pr(\lambda_1 \leq s_1, \ldots, \lambda_m \leq s_m)$$
(18)
$$= \sum_{\{i_1,\ldots,i_m\} \in I} \frac{1}{i_1! \cdots i_m!} \frac{\partial^{i_1 + \cdots + i_m}}{\partial z_1^{i_1} \cdots \partial z_m^{i_m}}$$
$$\times E \prod_{j=1}^{p}[1 + f_{z_1,\ldots,z_m}(\lambda_j)]\bigg|_{z_1 = \cdots = z_m = 0}.$$

By Proposition 1, we can equivalently say that $\Pr(\lambda_1 \leq s_1, \ldots, \lambda_m \leq s_m)$ can be expressed as the sum of a few coefficients in the power expansion of $\det(1 + K_{n,p} f_{z_1,\ldots,z_m})$.

Consider a rescaled and recentered version of the kernel of the operator $K_{n,p}$,

$$S_{n,p}(u, v) = \sigma_{n,p} K_{n,p}(\mu_{n,p} + \sigma_{n,p} u, \mu_{n,p} + \sigma_{n,p} v).$$

Let $S_{n,p}$ be an operator with kernel $S_{n,p}(u,v)$, which acts on $L^2((x_m,\infty))$. Note that

$$\det(1 + K_{n,p} f_{z_1,\ldots,z_m}) = \det(1 + S_{n,p} g_{z_1,\ldots,z_m}),$$

where $g_{z_1,\ldots,z_m}(x) = \sum_{k=1}^m (z_k - 1)\chi_{R_k}(x)$ and $R_1,\ldots,R_m$ equal $(x_1,\infty), (x_2,x_1], \ldots, (x_m, x_{m-1}]$, respectively.

Under the conditions of Proposition 2 and an additional condition that $n \geq p$, El Karoui [9] proves that there exists $\varepsilon > 0$ such that $S_{n,p}$ converges in trace-class norm to an operator $E \cdot Ai \cdot E^{-1}$ from the trace class, where $E$ is the operator which entails multiplication by $e^{-\varepsilon x}$ and $Ai$ is an integral operator acting on $L^2((x_m,\infty))$, which has the Airy kernel $Ai(x,y) = \int Ai(x+u) Ai(y+u)\,du$. The strategy of his proofs is the same as that of Baik, Ben Arous and Peche [4]. First, Proposition 1.2 of [4] is used to represent $S_{n,p}$ in the form $A_{n,p} \cdot B_{n,p}$, where $A_{n,p}$ and $B_{n,p}$ are operators with kernels $A_{n,p}(x,y) \equiv A_{n,p}(x+y) = \int_\Gamma e^{n f_{n,p}(z,x+y)}\,dz$ and $B_{n,p}(x,y) \equiv B_{n,p}(x+y) = \int_\Sigma e^{n g_{n,p}(z,x+y)}\,dz$. By the logic of steepest descent analysis, as $n$ and $p$ tend to infinity only the behavior of $f_{n,p}$ and $g_{n,p}$ around their respective maxima matters for the asymptotic behavior of $A_{n,p}(x)$ and $B_{n,p}(x)$. This fact is then exploited to show the convergence of $A_{n,p}(x)$ and $B_{n,p}(x)$ to $e^{-\varepsilon x} Ai(x)$ and $e^{\varepsilon x} Ai(x)$, respectively, which implies the convergence of $S_{n,p}$ to $E \cdot Ai \cdot E^{-1}$.

A careful reading of El Karoui's [9] proofs reveals that he needs the additional condition $n \geq p$ only to be able to use Proposition 1.2 of [4]. For all other purposes, the inequality $n/p \geq 1$ in his proofs can be replaced by $n/p > \gamma^2 > 0$ without changing the validity of the proofs. Therefore, El Karoui's [9] result and our Proposition 1 imply the convergence of $S_{n,p}$ to $E \cdot Ai \cdot E^{-1}$ without the extra condition that $n \geq p$.

Since trace-class operators form an ideal in the algebra of bounded linear operators, the operators $S_{n,p} \cdot g_{z_1,\ldots,z_m}$ and $E \cdot Ai \cdot E^{-1} \cdot g_{z_1,\ldots,z_m}$ must be from the trace class. Further, since $g_{z_1,\ldots,z_m}$ is a bounded function for all $z_1,\ldots,z_m$, $\|S_{n,p} \cdot g_{z_1,\ldots,z_m} - E \cdot Ai \cdot E^{-1} \cdot g_{z_1,\ldots,z_m}\|_1$ is less than or equal to $\|S_{n,p} - E \cdot Ai \cdot E^{-1}\|_1 \|g_{z_1,\ldots,z_m}\|$, which converges to zero as $n$ and $p$ tend to infinity. Here, $\|K\|_1$ denotes the trace-class norm of the operator $K$ and the above norm inequality follows from the inequalities of Theorem 1.6 in [22]. Hence, $S_{n,p} \cdot g_{z_1,\ldots,z_m}$ converges to $E \cdot Ai \cdot E^{-1} \cdot g_{z_1,\ldots,z_m}$ in the trace-class norm.

Now, since the Fredholm determinant is continuous with respect to the trace-class norm, $\det(1 + S_{n,p} \cdot g_{z_1,\ldots,z_m})$ converges to $\det(1 + E \cdot Ai \cdot E^{-1} \cdot g_{z_1,\ldots,z_m})$ for any $z_1,\ldots,z_m$. Further, since $E^{-1}$ and $g_{z_1,\ldots,z_m}$ commute, $\det(1 + E \cdot Ai \cdot E^{-1} \cdot g_{z_1,\ldots,z_m}) = \det(1 + E \cdot Ai \cdot g_{z_1,\ldots,z_m} \cdot E^{-1})$. However, determinants are invariant with respect to conjugation which leaves an operator in the trace class (see Remark 2.1 in [4]). Therefore, we have

$$\det(1 + S_{n,p} \cdot g_{z_1,\ldots,z_m}) \to \det(1 + Ai \cdot g_{z_1,\ldots,z_m})$$

for any $z_1, \ldots, z_m$.

Since $\det(1 + S_{n,p} \cdot g_{z_1,\ldots,z_m})$ exactly equals $E(\prod_{j=1}^{p}(1 + \sum_{k=1}^{m}(z_k - 1) \times \chi_{J_k}(\lambda_j)))$, it is a finite order polynomial in $z_1, \ldots, z_m$ and hence an analytic function of $z_1, \ldots, z_m$. Further, as follows, for example, from formulas (1.30) and (1.32) in [23],

$$(19) \quad |\det(1 + S_{n,p} \cdot g_{z_1,\ldots,z_m})| \leq \mathrm{const} \cdot \exp\left(\max_{j=1,\ldots,m} |z_j - 1| \|S_{n,p}\|_1\right)$$

(see also Lemma 3.3 of [22] for the case when $m = 1$). Since $S_{n,p}$ converges in the trace-class norm to $E \cdot Ai \cdot E^{-1}$, there exists a constant $M$ such that $\|S_{n,p}\|_1 < M$ for all $n$ and $p$. This fact, together with (19), implies that $\det(1 + S_{n,p} \cdot g_{z_1,\ldots,z_m})$ form a normal family of analytic functions (see [21], page 5). Hence, the convergence of $\det(1 + S_{n,p} \cdot g_{z_1,\ldots,z_m})$ to $\det(1 + Ai \cdot g_{z_1,\ldots,z_m})$ is uniform on any compact set in $\mathbb{C}^m$ and therefore all derivatives of $\det(1 + S_{n,p} \cdot g_{z_1,\ldots,z_m})$, and thus all derivatives of $E\prod_{j=1}^{p}[1 + f_{z_1,\ldots,z_m}(\lambda_j)]$, converge to the corresponding derivatives of $\det(1 + Ai \cdot g_{z_1,\ldots,z_m})$ at $z_1 = \cdots = z_m = 0$.

As shown by Johansson [12] [see his formulas (1.19), (3.46) and (3.48)], $F(x_1, \ldots, x_m)$, defined as

$$F(x_1, \ldots, x_m) = \sum_{\{i_1,\ldots,i_m\} \in I} \frac{1}{i_1! \cdots i_m!} \frac{\partial^{i_1 + \cdots + i_m}}{\partial z_1^{i_1} \cdots \partial z_m^{i_m}} \det(1 + Ai \cdot g_{z_1,\ldots,z_m})\bigg|_{z_1 = \cdots = z_m = 0},$$

is the distribution function for the $m$-dimensional joint Tracy–Widom distribution. Using (18) and the convergence of the derivatives of $E\prod_{j=1}^{p}[1 + f_{z_1,\ldots,z_m}(\lambda_j)]$ to those of $\det(1 + Ai \cdot g_{z_1,\ldots,z_m})$ just established, we conclude that

$$\Pr\left(\frac{\lambda_1 - \mu_{n,p}}{\sigma_{n,p}} \leq x_1, \ldots, \frac{\lambda_m - \mu_{n,p}}{\sigma_{n,p}} \leq x_m\right) \to F(x_1, \ldots, x_m).$$

Since $F(x_1, \ldots, x_m)$ is a continuous function, such a convergence implies that the joint distribution of $\sigma_{n,p}^{-1}(\lambda_1 - \mu_{n,p}), \ldots, \sigma_{n,p}^{-1}(\lambda_m - \mu_{n,p})$ weakly converges to the $m$-dimensional joint Tracy–Widom distribution. □

Concluding this section, we note that since [9] used the assumption that $n/p \geq 1$ only to be able to use Proposition 2.1 of [4], our finding that proposition holds for $n < p$ proves that all of El Karoui's [9] results are valid for $n/p$ remaining in a compact subset of $(0, \infty)$.

**4. Application.** In this section, we will see how Proposition 2 can be used in the analysis of excess stock return data generated by the approximate factor model (1). An approximate factor model for asset returns forms the core of Chamberlain and Rothschild's [8] extension of the arbitrage pricing theory (APT) of Ross [20]. The APT, one of the most important finance theories, shows that asset prices must be well explained by covariances of asset returns with a few common risk factors. An important practical question is how many such common factors exist.

This question has attracted considerable research attention. Roll and Ross ([19], page 1092) find that "at least three factors are important for pricing, but that it is unlikely that more than four are present." Brown and Weinstein ([5], page 713) "find evidence that there may be as few as 3 economywide factors, and certainly no more than 5 if the APM (arbitrage pricing model) is correct." Trzcinka [28] finds that there may be between one and five common risk factors. Connor and Korajczyk [7] report one or two factors in non-January months, but three to six factors for January returns. Huang and Jo ([10], page 988) find that " the evidence supports only a small number of factors, generally one and at most two." Bai and Ng [3] estimate the number of common factors in stock returns to be two. Two is also the preferred number of factors in [16]. Makarov and Papanikolaou [13] find evidence that there are four factors in stock returns.

In general, researchers find a small number of factors in the approximate factor model for excess stock returns. Often, they are uncertain about their point estimates. This uncertainty concerning point estimates is also reflected in the fact that different researchers often provide conflicting estimates. Even though the uncertainty is well recognized, it has never been formally quantified. Below, we will try to quantify this uncertainty. More precisely, we will find an asymptotic 95% confidence set for the number of factors by inverting a statistical test for the number of factors partially developed in a companion paper [17].

In the companion paper, I am interested in testing the null hypothesis of $k_0$ factors versus the alternative that the number of factors $k$ is larger than $k_0$ but smaller than $k_{\max}+1$, where $k_{\max}$ is an a priori maximum number of factors. I assume that the real-valued data $\{\xi_{it}\}_{1\leq i\leq p, 1\leq t\leq N}$, where $N=2n$, is generated by model (1), where the vectors of idiosyncratic components $\{\eta_{it}\}_{1\leq i\leq p}$ are i.i.d. Gaussian $N(0,\Sigma_p/2)$ and independent of factors $\{F_t\}_{1\leq t\leq N}$. To test the hypothesis, I propose first to construct a new complex-valued data set $\{\tilde{\xi}_{it}\}_{1\leq i\leq p, 1\leq t\leq n}$, where $\tilde{\xi}_{it} = \xi_{it} + \sqrt{-1}\xi_{i,t+n}$, and then to compute a test statistic $\max_{k_0<i\leq k_{\max}}(\gamma_i-\gamma_{i+1})/(\gamma_{i+1}-\gamma_{i+2})$, where $\gamma_i$ is the $i$th largest eigenvalue of the sample covariance matrix of the new dataset (this matrix is defined as $\frac{\tilde{X}\tilde{X}^*}{n/2}$, where $\tilde{X}=\{\tilde{\xi}_{it}\}_{1\leq i\leq p, 1\leq t\leq n}$). Theorem 3 in [17] shows that if the factors' explanatory power is strong enough so that $\gamma_1,\ldots,\gamma_k$ increase faster than $p^{2/3}$, then, under the null hypothesis, the asymptotic distribution of the proposed test statistic as $n$ and $p$ increase in such a way that $n/p$ remains in a compact subset of $(0,\infty)$ is the same as the asymptotic distribution of $\max_{0<i\leq k_{\max}-k_0}(\lambda_i - \lambda_{i+1})/(\lambda_{i+1} - \lambda_{i+2})$, where $\lambda_i$ is the $i$th largest eigenvalue of a $W_{\mathbb{C}}(\Sigma_p/n, n)$ matrix. Under the alternative, the test statistic explodes in probability as $n$ and $p$ increase.

Proposition 2 implies that the asymptotic distribution of $\max_{0<i\leq k_{\max}-k_0}(\lambda_i - \lambda_{i+1})/(\lambda_{i+1} - \lambda_{i+2})$ equals the distribution of $\max_{0<i\leq k_{\max}-k_0}(\mu_i - \mu_{i+1})/(\mu_{i+1} - \mu_{i+2})$, where $\mu_1, \ldots, \mu_{k_{\max}-k_0}$ have the joint $(k_{\max} - k_0)$-dimensional Tracy–Widom distribution. This result allows us to tabulate the critical values of the asymptotic distribution by using Monte Carlo simulations of large-dimensional matrices from the GUE to approximate the joint Tracy–Widom distribution. We approximate the joint 10-dimensional Tracy–Widom distribution of type two by the distribution of the 10 largest eigenvalues of a $1000 \times 1000$ matrix from the Gaussian Unitary Ensemble. We obtain an approximation for the latter distribution by simulating 30,000 independent matrices from the ensemble and numerically computing their first 10 eigenvalues.

Table 1 contains the critical values of our test for $k_{\max} - k_0 = 1, 2, \ldots, 8$. For example, in the table, the approximate 95% critical value of the test of three factors, versus the alternative $3 < k \leq 10$, is in the 5th row (counting from the bottom up) and the 2nd column (counting from the right). It equals 8.29.

Our test procedure can be interpreted as formalizing the widely used empirical method to determine the number of (classical) factors based on the visual inspection of the scree plot introduced by Cattell [6]. The scree plot is a line that connects the decreasing eigenvalues of the sample covariance matrix of the data plotted against their respective order numbers. In practice, it often happens that the scree plot shows a sharp break where the true number of factors ends and "debris" corresponding to the idiosyncratic influences appears. Our test statistic effectively measures the curvature of the scree plot at a would-be break point under the alternative hypothesis. When the alternative hypothesis is true, the curvature asymptotically tends to infinity. In contrast, under the null hypothesis, the curvature has a non-degenerate asymptotic distribution that does not depend on the model's parameter $\Sigma_p$.

TABLE 1
*Approximate percentiles of the test statistics for the tests of $k_0$ factors, versus the alternative of more than $k_0$ but less than $k_{\max} + 1$ factors*

| | $k_{\max} - k_0$ | | | | | | | |
|---|---|---|---|---|---|---|---|---|
| % | 1 | 2 | 3 | 4 | 5 | 6 | 7 | 8 |
| 50 | 1.27 | 1.95 | 2.30 | 2.54 | 2.74 | 2.92 | 3.09 | 3.24 |
| 60 | 1.53 | 2.24 | 2.59 | 2.88 | 3.10 | 3.31 | 3.49 | 3.65 |
| 70 | 1.86 | 2.61 | 3.01 | 3.32 | 3.59 | 3.82 | 4.01 | 4.20 |
| 80 | 2.37 | 3.19 | 3.65 | 4.02 | 4.32 | 4.59 | 4.83 | 5.05 |
| 85 | 2.75 | 3.62 | 4.15 | 4.54 | 4.89 | 5.20 | 5.45 | 5.70 |
| 90 | 3.33 | 4.31 | 4.91 | 5.40 | 5.77 | 6.13 | 6.42 | 6.66 |
| 91 | 3.50 | 4.49 | 5.13 | 5.62 | 6.03 | 6.39 | 6.67 | 6.92 |
| 92 | 3.69 | 4.72 | 5.37 | 5.91 | 6.31 | 6.68 | 6.95 | 7.25 |
| 93 | 3.92 | 4.99 | 5.66 | 6.24 | 6.62 | 7.00 | 7.32 | 7.59 |
| 94 | 4.20 | 5.31 | 6.03 | 6.57 | 7.00 | 7.41 | 7.74 | 8.04 |
| 95 | 4.52 | 5.73 | 6.46 | 7.01 | 7.50 | 7.95 | 8.29 | 8.59 |
| 96 | 5.02 | 6.26 | 6.97 | 7.63 | 8.16 | 8.61 | 9.06 | 9.36 |
| 97 | 5.62 | 6.91 | 7.79 | 8.48 | 9.06 | 9.64 | 10.11 | 10.44 |
| 98 | 6.55 | 8.15 | 9.06 | 9.93 | 10.47 | 11.27 | 11.75 | 12.13 |
| 99 | 8.74 | 10.52 | 11.67 | 12.56 | 13.42 | 14.26 | 14.88 | 15.25 |

To construct a 95% confidence interval for the number of common factors in excess stock returns, we use data provided by the Center for Research in Security Prices (CRSP) on monthly returns on $p = 972$ stocks traded on the NYSE, AMEX and NASDAQ during the period from January 1983 to December 2006. The data set includes those and only those companies for which CRSP provides monthly holding period return data for all months in the studied time interval. To obtain the excess returns on the stocks, we subtract the one-month risk-free rate provided by CRSP from the stock returns.

Since previous empirical research suggests that the number of common risk factors may be different in January and non-January months, we omit January data, which leaves $N = 264$ time observations of real-valued data. To obtain a complex-valued data set, we divide the real-valued data into two parts: the first containing all observations from February 1983 to December 1994 and the second containing all observations from February 1995 to December 2006. We then add the data from the first subperiod to the product of the imaginary unit and the data from the second subperiod. Hence, the dimensionality of the complex-valued data set is $p = 972$, $n = 132$. Note that $p$ is larger than $n$.

We maintain the assumption that the true number of factors is strictly less than seven. We take seven as an upper bound for the number of factors because it is consistent with conclusions of previous studies. If a 5%-asymptotic-size test of the null hypotheses of $j$ factors, versus the alternative that the number of factors is larger than $j$ but no larger than seven, does not reject the null hypothesis, we will include $j < 7$ in the asymptotic 95% confidence set for the number of factors.

Table 2 contains the first nine eigenvalues $\gamma_1, \ldots, \gamma_9$ of the sample covariance matrix of the complex-valued data set; the quantities $(\gamma_i - \gamma_{i+1})/(\gamma_{i+1} - \gamma_{i+2})$, $i = 1, \ldots, 7$; the test statistics $\max_{k_0 < i \leq 7}(\gamma_i - \gamma_{i+1})/(\gamma_{i+1} - \gamma_{i+2})$ for the tests of the null hypotheses of $k_0 = 0$, $k_0 = 1, \ldots, k_0 = 6$ factors, versus the alternatives that the number of factors is larger than $k_0$ but no larger than 7; and the corresponding 95% critical values taken from Table 1.

TABLE 2
*The largest eigenvalues of the sample covariance matrix of the complex-valued data and the test statistics for the tests of hypotheses that the number of factors equals $k_0$, versus the alternative of more than $k_0$ but no more than 7 factor*

| $i$ | 1 | 2 | 3 | 4 | 5 | 6 | 7 | 8 | 9 |
|---|---|---|---|---|---|---|---|---|---|
| $\gamma_i$ | 3.99 | 1.16 | 0.76 | 0.66 | 0.50 | 0.48 | 0.40 | 0.37 | 0.34 |
| $\frac{\gamma_i - \gamma_{i+1}}{\gamma_{i+1} - \gamma_{i+2}}$ | 7.14 | 3.77 | 0.65 | 12.73 | 0.15 | 2.98 | 1.18 | | |
| $k_0$ | 0 | 1 | 2 | 3 | 4 | 5 | 6 | | |
| $\max_{k_0 < i \leq 7} \frac{\gamma_i - \gamma_{i+1}}{\gamma_{i+1} - \gamma_{i+2}}$ | 12.73 | 12.73 | 12.73 | 12.73 | 2.98 | 2.98 | 1.18 | | |
| Critical values | 8.29 | 7.95 | 7.50 | 7.01 | 6.46 | 5.73 | 4.52 | | |

We reject the null hypotheses of 0, 1, 2 and 3 factors by the tests of asymptotic size 5%, but cannot reject the null hypotheses of 4, 5 and 6 factors. Hence, my 95% confidence set for the true number of common factors in the excess stock returns is $\{4, 5, 6\}$. This set intersects the range of estimates proposed by Brown and Weinstein [5], Roll and Ross [19] and Trzcinka [28]. It includes the point estimate 4 found in [13]. It is disjoint with the set of estimates reported by [7], [10] and [3]. The present author's previous estimate, 2, reported in [16], is not included in the set.

The 95% confidence set $\{4, 5, 6\}$ should, perhaps, appeal to the proponents of multifactor financial models. The good news is that 0 and 1 do not enter the set. However, this result should be taken with a grain of salt. The reason is that the rejection of the null hypothesis of 0, 1, 2 and 3 factors may be due to a failure of some of the primitive assumptions of the tests. For example, no one would truly believe that the idiosyncratic components of the excess stock returns are Gaussian and independent over time. Soshnikov's [24] results on the universality of the Tracy–Widom limit for the largest eigenvalue of sample covariance matrices require that $n/p$ approaches 1 as $n$ tends to infinity and assume that the tails of the distribution of the data points are relatively thin. For the financial data, the tails of the distribution may be heavy and the above test for the number of factors will then be invalid. In fact, the approximate factor model will no longer be a plausible description of the data because it assumes the existence of the second moments of the data.

Another discouraging possibility is that the asymptotics which the tests rely on may poorly approximate the finite-sample situation. Perhaps most importantly, the result of [17] that the asymptotic distribution of $\max_{k_0 < i \leq k_{\max}} (\gamma_i - \gamma_{i+1})/(\gamma_{i+1} - \gamma_{i+2})$ is the same as that of $\max_{0 < i \leq k_{\max} - k_0} (\lambda_i - \lambda_{i+1})/(\lambda_{i+1} - \lambda_{i+2})$, where $\lambda_i$ is the $i$th largest eigenvalue of a $W_{\mathbb{C}}(\Sigma_p/n, n)$ matrix, substantially uses the fact that if the true number of factors is $k_0$, then under the null hypothesis, $\gamma_{k_0}/\gamma_{k_0+1}$ increases faster than $p^{2/3}$. A casual inspection of the second row of Table 2 reveals that none of the ratios $\gamma_i/\gamma_{i+1}$ are large. Hence, to use the obtained 95% confidence set comfortably, one should check whether the asymptotic requirements in [17] can be relaxed. Such a check is left for future research.

**5. Conclusion.** In this paper, it has been shown that the joint distribution of the centered and normalized several largest eigenvalues of a $p$-dimensional complex Wishart matrix $W_{\mathbb{C}}(\Omega, n)$ converges to the joint Tracy–Widom distribution as $n$ and $p$ tend to infinity in such a way that $n/p$ remains in a compact subset of $(0, \infty)$. This result extends [4] and [9] in two directions. First, several of the largest eigenvalues, as opposed to the single largest eigenvalue, have been analyzed. Second, and most important, $n$ is allowed to be smaller than $p$, a situation corresponding to $W_{\mathbb{C}}(\Omega, n)$ being a singular matrix.

It has also been shown that all results of Baik, Ben Arous and Peche [4] and El Karoui [9] remain true if their assumption that $n/p$ remains in a compact subset of $[1, \infty)$ is replaced by a less restrictive assumption that $n/p$ remains in a compact subset of $(0, \infty)$.

Finally, it has been demonstrated how the theoretical result of this paper can be used to find a 95% confidence set for the number of common factors in excess stock returns. The established confidence set is $\{4, 5, 6\}$. Such a set formally quantifies the uncertainty concerning the true number of factors in excess stock returns evident in previous studies of the number of factors. The set supports some of the previous studies, but not others. Possible drawbacks in the proposed methodology of obtaining the 95% confidence set have been pointed out, suggesting some directions for future research.

## APPENDIX

PROOF OF (5). Consider the following analytic homeomorphism $g$ (see [11] for a useful summary of concepts from differential geometry) of almost all of the $U(p)$ on almost all the product $CV(p, n) \times U(p - n)$. Let $g(R) = \{Q_1, S\}$, where $Q_1 = R_1$, $S = H_{R_1}^* R_2$, $H_{R_1}$ is such that $[R_1, H_{R_1}] \in U(p)$ and the elements of $H_{R_1}$ are analytic functions of $R_1$. Columns of $H_{R_1}$ can, for example, be obtained by applying the Gram–Schmidt orthogonalization procedure to the projections of the first $p - n$ vectors of a fixed basis in $\mathbb{C}^p$ onto the $p - n$-dimensional subspace orthogonal to the columns of $R_1$. Such a construction will work for all $R$ such that the first $p - n$ vectors of the fixed basis and the columns of $R_1$ are linearly independent. Hence, it will work for almost all of the $U(p)$. The inverse of $g$ is given by $g^{-1}(\{Q_1, S\}) = [Q_1, H_{Q_1} S]$.

The homeomorphism $g$ maps the differential form $(R^* dR) \equiv (R_1^* dR_1) \times (R_2^* dR_1)(R_2^* dR_2)$ to the product of forms $(Q_1^* dQ_1)$, $(S^* H_{Q_1}^* dQ_1)$ and $(S^* H_{Q_1}^* d(H_{Q_1} S))$. We have $(Q_1^* dQ_1)(S^* H_{Q_1}^* dQ_1) = (Q_1^* dQ_1)|S| \times$ $(H_{Q_1}^* dQ_1) = (Q^* dQ_1)$, where $Q = [Q_1, H_{Q_1}]$. Further, the form $(S^* H_{Q_1}^* d(H_{Q_1} S))$ can be represented as a sum of $(S^* H_{Q_1}^* dH_{Q_1} S)$ and $(S^* dS)$. But the product of $(Q^* dQ_1)$ with $(S^* H_{Q_1}^* dH_{Q_1} S)$ is zero because $(Q^* dQ_1)$ is a form of maximum degree on the Stiefel manifold. Therefore,

$$\int_{R \in U(p)} e^{-\text{tr}(\Pi R_1 \Lambda R_1^*)} (R^* dR)$$
$$= \int_{S \in U(p-n)} \int_{Q_1 \in CV(p,n)} e^{-\text{tr}(\Pi Q_1 \Lambda Q_1^*)} (Q^* dQ_1)(S^* dS)$$
$$= \text{Vol}\{U(p-n)\} \int_{Q_1 \in CV(p,n)} e^{-\text{tr}(\Pi Q_1 \Lambda Q_1^*)} (Q^* dQ_1). \qquad \square$$

**A list of changes to the proofs in Baik, Ben Arous and Peche [4].** Our extension of Proposition 2.1 of [4] to the case where $n/p < 1$ ($M/N < 1$ in their notation) implies that the main results of that paper, contained in Theorems 1.1 and 1.2, hold if the requirement that $M/N$ remains in a compact subset of $[1, +\infty)$ as both $N$ and $M$ tend to infinity is relaxed to the requirement that $M/N$ remains in a compact subset of $(0, +\infty)$. Here is a list of few extra necessary changes to the proofs of Theorems 1.1 and 1.2. We use the notation of [4].

1. Formula (134) should be complemented by the following statement: Clearly, $T_1(t) > 0$ for $0 < \gamma < 1$.

2. The statement after formula (139) should be changed to $1 \leq \mu = (\gamma + 1)^2/\gamma^2 < \infty$.

3. Formula (141) should be changed to $\bar{\gamma} \leq \gamma \leq \gamma_0$ for fixed $\gamma_0 \geq 1$ and $0 < \bar{\gamma} \leq 1$.

4. Formula (143) should be changed to $0 < \delta < \min\{\frac{\nu^3}{6C_0}, \frac{\bar{\gamma}}{(1+\gamma_0)(1+\bar{\gamma})}\}$, $C_0 := \frac{(1+\gamma_0)^4(1+\bar{\gamma})^4}{4\gamma_0^4 \bar{\gamma}^4}(1 + \bar{\gamma}^2 \gamma_0^4)$.

5. Formula (144) should use definition of $C_0$ from change 4 above.

6. Formula (177) should be complemented by the following reasoning: Define

$$T_1(t) = (\gamma+1)^2 t^2 + (\gamma^2 - 1)t + 2\gamma.$$

Note that

$$\min_{t \in [0, 2p_c]} T_1(t) = \begin{cases} T_1\left(\frac{1-\gamma^2}{2(\gamma+1)^2}\right), & \text{for } 1/5 \leq \gamma \leq 1, \\ T_1(2p_c), & \text{for } 0 < \gamma < 1/5. \end{cases}$$

But $T_1(\frac{1-\gamma^2}{2(\gamma+1)^2}) \geq 6/25$ for $1/5 \leq \gamma \leq 1$ and $T_1(2p_c) = 6\gamma^2$ for $0 < \gamma < 1/5$. Hence, $T_1(t) > 0$ for $0 < \gamma < 1$.

7. Formula (237) should be changed in the same way as formula (141). See change 3 above.

8. Formula (239) should be changed to $0 < \delta < \min\{\frac{\Pi}{2}, \frac{1}{2(1+\gamma_0)}, \frac{\nu^2}{4C_1}\}$, $C_1 := \frac{8}{3}(\frac{1}{\Pi^3} + \frac{(1+\gamma_0)^3}{\bar{\gamma}^2})$.

9. Formula (241) should be changed to $|f^{(3)}(s)| = |\frac{2}{s^3} - \frac{2}{\gamma^2(s-1)^3}| \leq \frac{2}{(\pi_1-\delta)^3} + \frac{2}{\bar{\gamma}^2(1-\pi_1-\delta)^3} \leq \frac{16}{\Pi^3} + \frac{16(1+\gamma_0)^3}{\bar{\gamma}^2} = 6C_1$.

**Acknowledgments.** I am very grateful to Jinho Baik for his encouraging interest and useful comments. I am thankful to an anonymous referee for valuable suggestions on how to improve the presentation of my results.

## REFERENCES

[1] ANDERSON, T. W. (1984). *An Introduction to Multivariate Statistical Analysis.* Wiley, New York. MR0771294
[2] ANDREIEF, C. (1883). Note sur une relation les integrales definies des produits des fonctions. *Mem. de la Soc. Sci. Bordeaux* **2** 1–14.
[3] BAI, J. and NG, S. (2002). Determining the number of factors in approximate factor models *Econometrica* **70** 191–221. MR1926259
[4] BAIK, J., BEN AROUS, G. and PECHE, S. (2005). Phase transition of the largest eigenvalue for nonnull complex sample covariance matrices. *Ann. Probab.* **33** 1643–1697. MR2165575
[5] BROWN, S. and WEINSTEIN, M. (1983). A new approach to testing asset pricing models: The bilinear paradigm. *J. Finance* **38** 711–743.
[6] CATTELL, R. B. (1966). The scree test for the number of factors. *Multivariate Behavioral Research* **1** 245–276.
[7] CONNOR, G. and KORAJCZYK, R. (1993). A test for the number of factors in an approximate factor model. *J. Finance* **58** 1263–1291.
[8] CHAMBERLAIN, G. and ROTHSCHILD, M. (1983). Arbitrage, factor structure, and mean-variance analysis on large asset markets. *Econometrica* **51** 1281–1304. MR0736050
[9] EL KAROUI, N. (2007). Tracy–Widom limit for the largest eigenvalue of a large class of complex sample covariance matrices. *Ann. Probab.* **35** 663–714. MR2308592
[10] HUANG, R. and JO, H. (1995). Data frequency and the number of factors in stock returns. *J. Banking and Finance* **19** 987–1003.
[11] JAMES, A. (1954). Normal multivariate analysis and the orthogonal group. *Ann. Math. Statist.* **25** 40–75. MR0060779
[12] JOHANSSON, K. (2001). Discrete orthogonal polynomial ensembles and the Plancherel measure. *Ann. of Math.* **153** 259–296. MR1826414
[13] MAKAROV, I. and PAPANIKOLAOU, D. (2007). Sources of systematic risk. Unpublished manuscript, London Business School.
[14] MEHTA, M. (2004). *Random Matrices.* Academic Press, New York. MR2129906
[15] OLVER, F. (1974). *Asymptotics and Special Functions.* Academic Press, New York. MR0435697
[16] ONATSKI, A. (2005). Determining the number of factors from the empirical distribution of eigenvalues. Manuscript, Columbia Univ.
[17] ONATSKI, A. (2007). A formal statistical test for the number of factors in the approximate factor models. Unpublished manuscript, Columbia Univ.
[18] RATNARAJAH, T. and VAILLANOURT, R. (2005). Complex singular Wishart matrices and applications. *Comput. Math. Appl.* **50** 399–411. MR2165429

[19] ROLL, R. and ROSS, S. (1980). An empirical investigation of the arbitrage pricing theory. *J. Finance* **5** 1073–1103.
[20] ROSS, S. (1976). The arbitrage theory of capital asset pricing. *J. Econom. Theory* **13** 341–360. MR0429063
[21] RUDIN, W. (1980). *Function Theory in the Unit Ball of* $\mathbf{C}^n$. Springer, New York. MR0601594
[22] SIMON, B. (2005). *Trace Ideals and Their Applications*, 2nd ed. Amer. Math. Soc., Providence, RI. MR2154153
[23] SOSHNIKOV, A. (2000). Determinantal random point fields. *Russian Math. Surveys* **55** 923–975. MR1799012
[24] SOSHNIKOV, A. (2002). A note on universality of the distribution of the largest eigenvalues in certain sample covariance matrices. *J. Statist. Phys.* **108** 1033–1056. MR1933444
[25] STOCK, J. and WATSON, M. (2002). Forecasting using principal components from a large number of predictors. *J. Amer. Statist. Assoc.* **97** 1167–1179. MR1951271
[26] TRACY, C. A. and WIDOM, H. (1994). Level spacing distributions and the Airy kernel. *Comm. Math. Phys.* **159** 151–174. MR1257246
[27] TRACY, C. A. and WIDOM, H. (1998). Correlation functions, cluster functions, and spacing distributions for random matrices. *J. Statist. Phys.* **92** 809–835. MR1657844
[28] TRZCINKA, C. (1986). On the number of factors in the arbitrage pricing model. *J. Finance* **41** 347–368.

DEPARTMENT OF ECONOMICS
COLUMBIA UNIVERSITY
NEW YORK, NEW YORK 10027
USA
E-MAIL: ao2027@columbia.edu